\documentclass[11pt]{article}
\usepackage{amsfonts}
\usepackage{mathrsfs}
\usepackage{amsmath}
\usepackage{amssymb}
\usepackage{arydshln}
\usepackage{graphicx}
\usepackage{epsfig}
\usepackage{epic}
\usepackage{cite}
\usepackage{setspace}
\usepackage{amsthm,latexsym}
\usepackage{float}
\usepackage[utf8]{inputenc}

\renewcommand{\paragraph}{\roman{paragraph}}
\newtheorem{theorem}{ \bf Theorem}[section]
\newtheorem{lemma}[theorem]{ \bf Lemma}

\newtheorem{corol}[theorem]{ \bf Corollary}

\begin{document}

\title{\sf On the multiplicity of 1 as a Laplacian eigenvalue of a graph}
\author{ \ Fenglei Tian$^{1,}$\thanks{Corresponding author. E-mail address: tflqsd@qfnu.edu.cn.  Supported by ``the National Natural Science Foundation of China (No.12101354, 12371025), the Natural Science Foundation of Shandong Province (No. ZR2024MA032) and the Youth Innovation Team Project of Shandong Province Universities (No. 2023KJ353)''. }\ ,\ \
Dein Wong$^2$ \ \
~~\\
\noindent{\small\it \ 1.\  School of Management, Qufu Normal University, Rizhao, China.}\\
\noindent{\small\it \ 2.\ School of Mathematics, China University of Mining and Technology, Xuzhou, China.}
}
\date{}
\maketitle
\noindent {\bf Abstract:} \ Let $G$ be a graph with $p(G)$ pendant vertices and $q(G)$ quasi-pendant vertices. Denote by $m_{L(G)}(\lambda)$ the multiplicity of $\lambda$ as a Laplacian eigenvalue of $G$. Let $\overline{G}$ be the reduced graph of $G$, which can be obtained from $G$ by deleting some pendant vertices such that $p(\overline{G})=q(\overline{G})$. We first prove that $m_{L(G)}(1)=p(G)-q(G)+m_{L(\overline{G})}(1)$.
Since deleting pendant path $P_3$ does not change the multiplicity of Laplacian eigenvalue 1 of a graph, we further  focus on reduced graphs without pendant path $P_3$. Let $T$ be a reduced tree on $n(\geq 6)$ vertices without pendant path $P_3$, then it is proved that
$$m_{L(T)}(1)\leq \frac{n-6}{4},$$
and all the  trees attaining the upper bound are characterized completely. As an application, for a reduced unicyclic graph $G$ of order $n\geq 10$ without pendant path $P_3$, we get
$$m_{L(G)}(1)\leq \frac{n}{4},$$
and all the unicyclic graphs attaining the upper bound are determined completely.

\section{Introduction}

\quad All graphs considered in this paper are simple and undirected. Let $G$ be a graph with vertex set $V(G)$ and edge set $E(G)$. The line graph of $G$ is denoted by $G^l$, whose vertex set $V(G^l)$ is the edge set $E(G)$ of $G$. If two edges in $G$ are incident, then the corresponding two vertices in $G^l$ are adjacent.  The degree $d_u$ of a vertex $u$ is the number of vertices in the neighbor set $N_G(u)$. A vertex $u$ is called a pendant vertex if $d_u=1$, and the vertex adjacent to $u$ is called  a  quasi-pendant vertex. The number of pendant vertices (resp., quasi-pendant vertices) of $G$ is denoted by $p(G)$ (resp., $q(G)$). Denote by $|G|$ the order of $G$ for brevity. We write $u\thicksim v$ if $u$ and $v$ are adjacent. The distance of $u$ and $v$ is denoted by $d(u,v)$ and the diameter of a graph $G$ is written as $diam(G)$. Denote by $I$ the identity matrix.
For a subset $S\subset V(G)$, let $G-S$ be the induced subgraph obtained from $G$ by deleting the vertices of $S$ with the incident edges. Sometimes, we write $G-G_1$ instead of $G-V(G_1)$ for convenience, if $G_1$ is an induced subgraph of $G$. Denote by  $G-e$ a spanning subgraph of $G$ obtained by removing the edge $e$ from $G$. A graph $G$ is called a {\it reduced graph}, if $p(G)=q(G)$ in $G$. For a given graph $G$, one can obtain the reduced graph of $G$ by deleting some pendant vertices until each quasi-pendant vertex possesses exactly one pendant vertex.
The multiplicity of an eigenvalue $\lambda$ of a matrix $M$ is denoted by $m_M(\lambda)$.
The Laplacian matrix $L(G)=D(G)-A(G)$ of a graph $G$ is a well-known matrix in spectral graph theory, where $A(G)$ is the adjacency matrix of $G$ and $D(G)$ is the degree diagonal matrix of $G$. The eigenvalues of $L(G)$ are called the Laplacian eigenvalues of $G$.

There have been so many papers about the Laplacian eigenvalues of graphs. Recently, the multiplicity of Laplacian eigenvalues of a graph has attracted much attention (see \cite{Akbari,Andrade,Barik,Faria,Guo,Gutman,HuangQ} for example). Here we concentrate on the multiplicity of the Laplacian eigenvalue 1 of graphs.
Faria\cite{Faria} showed that for a graph $G$,
$m_{L(G)}(1)\geq p(G)-q(G).$ The authors of \cite{Andrade} proved that the above bound is attained if every vertex of degree at least two is a quasi-pendant vertex in $G$.
If $\lambda >1$ is an integral Laplacian eigenvalue of a tree $T$ on $n$ vertices, Grone, Merris and Sunder \cite{Merris1} proved that $\lambda$ divides $n$ and $m_{L(T)}(\lambda)=1$. For an arbitrary Laplacian eigenvalue $\lambda$ of a tree $T$, they \cite{Merris1} showed that
\begin{equation}\label{e1}
m_{L(T)}(\lambda)\leq p(T)-1.
\end{equation}
All the trees attaining the upper bound in (\ref{e1}) were  characterized completely by Yang and Wang \cite{WangL}.
Let $N$ be the submatrix of $L(G)$ indexed by the vertices which are neither pendant vertices nor quasi-pendant vertices of a graph $G$, then the authors of \cite{Merris1,Andrade} proved that
\begin{equation}\label{q2}
m_{L(G)}(1)=p(G)-q(G)+m_N(1).
\end{equation}
Guo, Feng and Zhang \cite{Guo} considered the effect of adding an edge between two disjoint graphs on the multiplicity of Laplacian eigenvalues. Applying the effect they characterized all trees of order $n$ with $n-6 \leq m_{L(T)}(1)\leq n$.
Let $T$ be a tree obtained by adding an edge between a vertex of a tree $T_1$ and a vertex of a tree $T_2$.  Barik, Lal and Pati\cite{Barik} discussed the relationship among $m_{L(T)}(1)$, $m_{L(T_1)}(1)$ and $m_{L(T_2)}(1)$. Wen and Huang \cite{HuangQ} investigated the multiplicity of Laplacian eigenvalue 1 of unicyclic graphs on $n$ vertices.

Moreover, if $G$ is obtained from a graph $G_1$ and a path $P_k$ by joining a pendant vertex of $P_k$ and an arbitrary vertex of $G_1$, then we say $P_k$ is a pendant path of $G$. From \cite{Merris1}, we know that deleting a pendant path $P_3$ does not change  the multiplicity of Laplacian eigenvalue 1 of a graph.  In \cite{Tian}, the authors proved that $m_{L(T)}(1)=p(T)-q(T)+m_{L(\overline{T})}(1),$ where $\overline{T}$ is the reduced tree of $T$.
Inspired by these conclusions, to investigate the multiplicity of Laplacian eigenvalue 1 of trees, we restrict our consideration to the trees with $p(T)=q(T)$ and containing no pendant path $P_3$.
The main results of this paper are presented below.

\begin{theorem}\label{mainth0}\ Let $G$ be a graph of order $n(\geq 4)$. Denote by $\overline{G}$ the reduced graph of $G$. Then it is obtained that
\begin{equation}\label{q3}
m_{L(G)}(1)=p(G)-q(G)+m_{L(\overline{G})}(1).
\end{equation}
\end{theorem}

By Theorem \ref{mainth0}, to investigate the multiplicity of Laplacian eigenvalue 1 of a graph, we could restrict our attention  to the reduced graphs. Compared with (\ref{q3}), the above matrix $N$ in (\ref{q2}) is not a Laplacian matrix of some graph, so we have very few combinatorial methods to deal with $m_N(1)$. Let $\mathcal{G}(n)$ be the set of reduced graphs of order $n$ without pendant path $P_3$.

\begin{theorem}\label{mainth1}\ Let $T\in \mathcal{G}(n)$ be a tree on $n (\geq 6)$ vertices. Then we obtain that $$m_{L(T)}(1)\leq \frac{n-6}{4},$$
and the equality holds if and only if $T$ is isomorphic to the tree in Fig. 1.
\end{theorem}

\begin{figure}[htbp]
  \centering
  \setlength{\abovecaptionskip}{0.5cm} 
  \setlength{\belowcaptionskip}{0pt}
  \includegraphics[width=3.7 in]{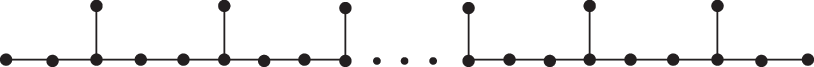}
  \caption{The extremal tree in Theorem \ref{mainth1}.}
\end{figure}

\begin{theorem}\label{mainth2}\ Let $G\in \mathcal{G}(n)$ be a unicyclic graph on $n(\geq 10)$ vertices.  Then we have
$$m_{L(G)}(1)\leq \frac{n}{4},$$
and the equality holds if and only if $G$ is isomorphic to the graph in Fig. 2.
\end{theorem}

\begin{figure}[htbp]
  \centering
  \setlength{\abovecaptionskip}{0.5cm} 
  \setlength{\belowcaptionskip}{0pt}
  \includegraphics[width=1 in]{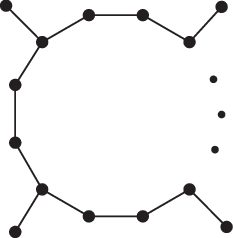}
  \caption{The extremal graph in Theorem \ref{mainth2}.}
\end{figure}

The rest of the paper is arranged as follows. In Section 2, some lemmas are introduced. In Section 3, we will show the proofs of Theorems \ref{mainth0}, \ref{mainth1} and \ref{mainth2}.

\section{ Preliminaries }
\quad
The following lemma gives us another indirect angle to investigate the multiplicity of Laplacain eigenvalue 1 of a bipartite graph.
\begin{lemma}{\rm\cite{Tian}}\label{eqlemma}\ Let $G$ be a bipartite graph. Then the multiplicity of $1$ as a Laplacian eigenvalue of $G$ is equal to the multiplicity of $-1$ as an adjacency eigenvalue of $G^l$, that is,$$m_{L(G)}(1)=m_{A(G^l)}(-1).$$ \end{lemma}

\begin{lemma}{\rm\cite{Merris1}}\label{boundlemma}\ Let $G-e$ be the graph obtained from a graph $G$ by deleting an edge $e$. Then
$$-1+m_{L(G-e)}(1)\leq m_{L(G)}(1)\leq 1+ m_{L(G-e)}(1).$$
\end{lemma}

\begin{lemma}{\rm\cite{Merris1}}\label{path3} Let $T$ be a tree obtained from a tree $T_1$ and a path $P_3$ by joining a vertex of $T_1$ to a pendant vertex of $P_3$. Then it follows that
$$m_{L(T)}(1)=m_{L(T_1)}(1).$$
\end{lemma}

\begin{lemma}\label{mainlemma} Let $G$ be a graph with a pendant vertex $u$ and $u\thicksim v$. Suppose that $d_v\geq 3$ and $w\thicksim v$, then let $\Gamma$ be the graph obtained from $G-e_{vw}$ and a disjoint path $P_2=x\thicksim y$ by joining $w$ and $y$. Then we have
$$m_{L(G)}(1)=m_{L(\Gamma)}(1).$$
\end{lemma}

\vskip 2mm
\noindent
{\bf Proof.} By arranging the vertices of $G$, we could write $L(G)-I$ as below such that the leading principal submatrix is indexed by $\{u, v, w\}$,
$$L(G)-I=
\left(
  \begin{array}{cccc}
    0 & -1 & 0 & \bf{0} \\
    -1 & d_v-1 & -1 & \alpha^T \\
    0 & -1 & d_w-1 & \beta^T \\
    \bf{0} & \alpha & \beta & B \\
  \end{array}
\right).$$
Set an invertible matrix $P$ such that
$$P=\left(
    \begin{array}{cccc}
      1 & \frac{1}{2} & -1 & \bf{0} \\
      0 & 1 & 0 & \bf{0} \\
      0 & 0 & 1 & \bf{0} \\
      \bf{0} & \bf{0} & \bf{0} & I_{n-3} \\
    \end{array}
  \right),$$
then it follows that
$$P^T(L(G)-I)P=
\left(
  \begin{array}{cccc}
    0 & -1 & 0 & \bf{0} \\
    -1 & d_v-2 & 0 & \alpha^T \\
    0 & 0 & d_w-1 & \beta^T \\
    \bf{0} & \alpha & \beta & B \\
  \end{array}
\right).$$
Write $L(\Gamma)-I$ as below such that the leading principal submatrix is indexed by $\{x, y, u, v, w\}$,
$$L(\Gamma)-I=
\left(
  \begin{array}{cc:cccc}
    0 & -1 & 0 & 0 & 0 & \bf{0} \\
    -1 & 1 & 0 & 0 & -1 & \bf{0} \\
    \hdashline
    0 & 0 & 0 & -1 & 0 & \bf{0} \\
    0 & 0 & -1 & d_v-2 & 0 & \alpha^T \\
    0 & -1 & 0 & 0 & d_w-1 & \beta^T \\
    \bf{0} & \bf{0} & \bf{0} & \alpha & \beta & B \\
  \end{array}
\right)$$
It is not hard to see that $P^T(L(G)-I)P$ is a principal submatrix of $L(\Gamma)-I$ and
$$m_{L(\Gamma)-I}(0)=m_{P^T(L(G)-I)P}(0).$$
Since $rank(P^T(L(G)-I)P)=rank(L(G)-I)$, then
$$m_{P^T(L(G)-I)P}(0)=m_{L(G)-I}(0).$$
Hence, it follows that $m_{L(\Gamma)-I}(0)=m_{L(G)-I}(0)$, which implies that $m_{L(\Gamma)}(1)=m_{L(G)}(1)$, as required.
\hfill$\square$

\vskip 2mm
The following technique of reduction operation for a tree has been given in \cite{Tian}. Here we generalize it to an arbitrary graph.

\noindent
{$\bullet$ {\bf Reduction Operation }}\ \ Let $u$ be a pendant vertex of a graph $G$ and $u\thicksim v$. Let $N_G(v)\setminus \{u\}=\{v_1, v_2, \cdots, v_s\}$ be the neighbors, except $u$, of $v$. For graph $G-\{u, v\}$, joining each $v_i\ (1\leq i\leq s)$ to an endpoint of a path $P_2$, the resultant graph is denoted by $G'$ and is called the {\it reduction graph} of $G$. An example of doing reduction operation to a graph $G$ is illustrated in Fig. 3.

It is not hard to see that if $d_v=2$ in $G$, then $G'=G$. One can do reduction operation to $G'$ again, until $G'$ contains no quasi-pendant vertex of degree greater than 2. At this time, we call $G'$ the {\it final reduction graph} of $G$.

\begin{figure}[htbp]
  \centering
  \setlength{\abovecaptionskip}{0.5cm} 
  \setlength{\belowcaptionskip}{0pt}
  \includegraphics[width=5 in]{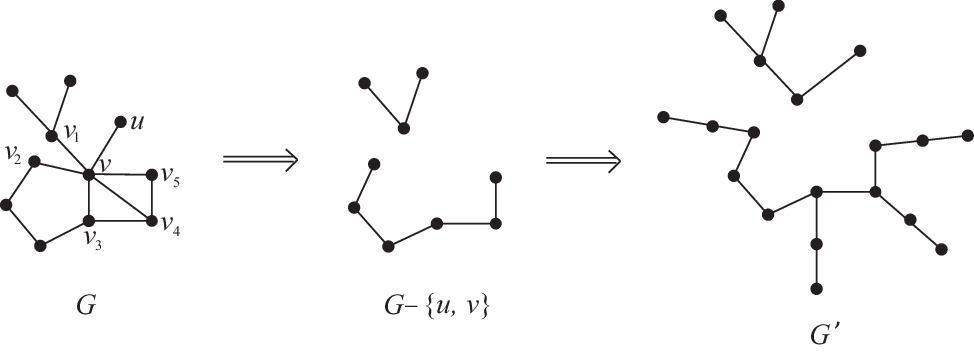}
  \caption{An example of doing reduction operation to a graph.}
\end{figure}

Let $v$ be a quasi-pendant vertex of a graph $G$ with $d_v\geq 3$. Then by using Lemma \ref{mainlemma} recursively until $d_v=2$, the following corollary is clear.

\begin{corol}\label{corol}\ Let $u$ be a pendant vertex of a graph $G$ and $u\thicksim v$ with $d_v\geq 3$. Let $G'$ be the reduction graph of $G$ with respect to $u$ and $v$. Then
$$m_{L(G)}(1)= m_{L(G')}(1).$$
\end{corol}



\begin{lemma}{\rm\cite{Tian}}\label{previousth}\ Let $T$ be a reduced tree on $n\geq 6$ vertices. Let $T_i'$ $(i=1, 2,\cdots, k)$ be the components of the final reduction graph of $T$. Then
$$m_{L(T)}(1)\leq \frac{n-2}{4},$$
and the equality holds if and only if the following two assertions hold:\\
(i)\  Each vertex of degree greater than $2$ is a quasi-pendant vertex in $T$;\\
(ii)\ Each $T_i'$ of $T$ is a path $P_6$.\end{lemma}

Let $P_k=u_1u_2\cdots u_k$ be an {\it internal path} of a graph $G$, that is, $u_1\thicksim u_2\thicksim\cdots \thicksim u_k$ and  $d_{u_i}=2$ ($2\leq i\leq k-1$), and $u_1, u_k$ have no common neighbors in $G$. The following result was presented in \cite{WongD}, but there is a flaw in their proof. Here we provide an amended proof.

\begin{lemma}{\rm\cite{WongD}}\label{innerpath}\ Let $P_4=u_1u_2u_3u_4$ be an internal path of a graph $G$. Let $H$ be the graph obtained from $G$ by contracting $P_4$ to a vertex $v$. Then, for the multiplicity of $-1$ as an eigenvalue of the  adjacency matrix $A(G)$, we have $m_{A(G)}(-1)=m_{A(H)}(-1)$.
\end{lemma}

\vskip 2mm
\noindent
{\bf Proof.} By arranging the vertices of $G$, we could write $A(G)+I$ as below such that the leading principal submatrix is indexed by $\{u_1, u_2, u_3, u_4\}$,
$$A(G)+I=\left(
    \begin{array}{ccccc}
      1 & 1 & 0 & 0 & \alpha^T \\
      1 & 1 & 1 & 0 & {\bf 0} \\
      0 & 1 & 1 & 1 & {\bf 0} \\
      0 & 0 & 1 & 1 &  \beta^T \\
      \alpha & {\bf 0} & {\bf 0} &  \beta & A(G-P_4)+I \\
    \end{array}
  \right).$$
By the technique of elementary transformation to $A(G)+I$, we set $$P=\left(
  \begin{array}{ccccc}
    1 & -1 & 0 & 0 & {\bf 0} \\
    0 & 1 & 0 & 0 & {\bf 0} \\
    0 & -1 & 1 & 0 & {\bf 0} \\
    1 & -1 & 0 & 1 & {\bf 0} \\
    {\bf 0} & -\alpha & \alpha & {\bf 0} & I \\
  \end{array}
\right),$$
and then
$$P(A(G)+I)P^T=\left(
                 \begin{array}{ccccc}
                   0 & 0 & -1 & 0 & {\bf 0} \\
                   0 & 1 & 0 & 0 & {\bf 0} \\
                   -1 & 0 & 0 & 0 & {\bf 0} \\
                   0 & 0 & 0 & 1 & \alpha^T+\beta^T \\
                   {\bf 0} & {\bf 0} & {\bf 0} & \alpha+\beta &  A(G-P_4)+I \\
                 \end{array}
               \right).$$
It is not hard to see that
$$\left(
    \begin{array}{cc}
      1 & \alpha^T+\beta^T \\
      \alpha+\beta &  A(G-P_4)+I \\
    \end{array}
  \right)=A(H)+I.
$$
Thus $rank(A(G)+I)=rank(P(A(G)+I)P^T)=3+rank(A(H)+I)$, which implies that $m_{A(G)}(-1)=m_{A(H)}(-1)$.
\hfill$\square$

\vskip 2mm
If the graph $G$ in Lemma \ref{innerpath} is the line graph of a tree $T$, then the following corollary is clear.

\begin{corol}\label{innerpathcorol}\ Let $T$ be a tree with an internal path $P_5=u_1u_2u_3u_4u_5$. Let $H$ be the tree obtained from $T$ by deleing the vertices $\{u_2, u_3, u_4\}$ and joining $u_1$ and $u_5$ with an edge. Then it is obtained that $m_{L(T)}(1)=m_{L(H)}(1)$.
\end{corol}

\begin{lemma}\label{extremallemma}\ Let $T$ be a tree on $n$ vertices as depicted in Fig. $4$, which has a local structure induced by $\{u_1, u_2, \cdots, u_8\}$. Let $T_1$ be the tree obtained from $T$ by deleting $\{u_1, u_2, u_3, u_4\}$ (see Fig. $4$). Then we get $m_{L(T)}(1)=1 + m_{L(T_1)}(1)$.
\end{lemma}

\begin{figure}[htbp]
  \centering
  \setlength{\abovecaptionskip}{0cm}
  \setlength{\belowcaptionskip}{0pt}
  \includegraphics[width=4.3 in]{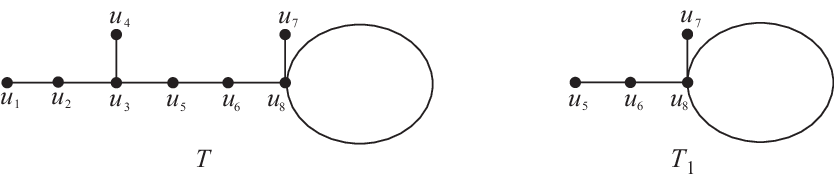}
  \caption{The trees in Lemma \ref{extremallemma}.}
\end{figure}

\vskip 2mm
\noindent
{\bf Proof.} \ By Lemma \ref{eqlemma}, we consider the adjacency matrix of the line graph of $T$ instead of the Laplacian matrix of $T$. For brevity, denote $\Gamma = T-\{u_1, u_2, \cdots, u_7\}$.  The matrix $A(T^l)+I$ can be written as
$$A(T^l)+I=\left(
             \begin{array}{ccccccccc}
               1 & 1 & 0 & 0 & 0 & 0 & 0  & {\bf 0} \\
               1 & 1 & 1 & 1 & 0 & 0 & 0  & {\bf 0} \\
               0 & 1 & 1 & 1 & 0 & 0 & 0  & {\bf 0} \\
               0 & 1 & 1 & 1 & 1 & 0 & 0  & {\bf 0} \\
               0 & 0 & 0 & 1 & 1 & 1 & 0  & {\bf 0} \\
               0 & 0 & 0 & 0 & 1 & 1 & 1  & \alpha^T \\
               0 & 0 & 0 & 0 & 0 & 1 & 1  & \alpha^T \\
               {\bf 0} & {\bf 0} & {\bf 0} & {\bf 0} & {\bf 0} & \alpha & \alpha & A(\Gamma^l) \\
             \end{array}
           \right).$$
Set the matrix $P$ as below and calculate $P^T(A(T^l)+I)P$,

$$P=\left(
      \begin{array}{cccccc}
        I_3 & {\bf 0} & {\bf 0} & {\bf 0} & {\bf 0} & {\bf 0} \\
        {\bf 0} & 1 & 0 & 0 & 0 & {\bf 0} \\
        {\bf 0} & 0 & 1 & 0 & 0 & {\bf 0} \\
        {\bf 0} & -1 & 0 & 1 & 0 & {\bf 0} \\
        {\bf 0} & 1 & 0 & -1 & 1 & {\bf 0} \\
        {\bf 0} & {\bf 0} & {\bf 0} & {\bf 0} & {\bf 0} & I_{n-8} \\
      \end{array}
    \right)\ \ \text{and}\ \
    P^T(A(T^l)+I)P=\left(
                     \begin{array}{cc}
                       A_1 & {\bf 0} \\
                       {\bf 0} & A_2 \\
                     \end{array}
                   \right),$$
where $A_1=\left(
             \begin{array}{cccc}
               1 & 1 & 0 & 0 \\
               1 & 1 & 1 & 1 \\
               0 & 1 & 1 & 1 \\
               0 & 1 & 1 & 1 \\
             \end{array}
           \right)$ with $rank(A_1)=3$ obviously and
$A_2=\left(
       \begin{array}{cccc}
         1 & 1 & 0  & {\bf 0} \\
         1 & 0 & 0  & {\bf 0} \\
         0 & 0 & 1  & \alpha^T \\
         {\bf 0} & {\bf 0} & \alpha & A(\Gamma^l) \\
       \end{array}
     \right)$.
It is clear that $A_2$ is congruent to the adjacency matrix $A(T_1^l)$ of the line graph $T_1^l$ of $T_1$, and so $rank(A_2)=rank(A(T_1^l))$. As a result,
$$rank(A(T^l)+I)=rank(P^T(A(T^l)+I)P)=rank(A_1)+rank(A_2)=3+rank(A(T_1^l)),$$
which indicates that $m_{A(T^l)}(-1)=1 + m_{A(T_1^l)}(-1)$. By Lemma \ref{eqlemma}, $m_{L(T)}(1)=1 + m_{L(T_1)}(1)$ follows immediately.\hfill$\square$

\vskip 2mm
For a tree $T$ on $n\geq 5$ vertices, if there exists a center vertex (say $u$) and each component of $T-u$ is $P_2$, then we call $T$ a {\it star-like tree}.

\begin{lemma}\label{starlike}\ Let $T$ be a star-like tree. Let $H$ be a tree obtained from two star-like trees by joining their center vertices with an edge. Then we have\\
(i){\rm\cite{Guo,Tian}} \ $m_{L(T)}(1)=0$;\\
(ii) \ $m_{L(H)}(1)=0$.
\end{lemma}

\vskip 2mm
\noindent
{\bf Proof.} \ The conclusion of (i) can be seen in \cite{Guo,Tian}. Here we show the proof of (ii). Let $T_1$ and $T_2$ be two star-like trees on $2s+1$ ($s\geq 2$) and $2t+1$ ($t\geq 2$) vertices, respectively. Let $H$ be obtained from $T_1$ and $T_2$ by joining their center vertices with an edge. Considering the adjacency matrix of the line graph $H^l$ of $H$, we could write
$$A(H^l)+I=\left(
    \begin{array}{ccccc}
      I_s & I_s & {\bf 0} & {\bf 0} & {\bf 0} \\
      I_s & J_s & {\bf 1} & {\bf 0} & {\bf 0} \\
      {\bf 0} & {\bf 1}^T & 1 & {\bf 1}^T & {\bf 0} \\
      {\bf 0} & {\bf 0} & {\bf 1} & J_t & I_t \\
      {\bf 0} & {\bf 0} & {\bf 0} & I_t & I_t \\
    \end{array}
  \right),$$
where $J_k$ is the all-ones matrix of order $k$ and ${\bf 1}$ is an  all-ones column vector. Set an invertible matrix $P$ such that
$$P=\left(
    \begin{array}{ccccc}
      I_s & -I_s & \frac{1}{s-1}{\bf 1} & {\bf 0} & {\bf 0} \\
      {\bf 0} & I_s & \frac{-1}{s-1}{\bf 1} & {\bf 0} & {\bf 0} \\
      {\bf 0} & {\bf 0} & 1 & {\bf 0} & {\bf 0} \\
      {\bf 0} & {\bf 0} & \frac{-1}{t-1}{\bf 1} & I_t & {\bf 0} \\
      {\bf 0} & {\bf 0} & \frac{1}{t-1}{\bf 1} & -I_t & I_t \\
    \end{array}
  \right),$$
then we obtain that
$$P^T(A(H^l)+I)P=
\left(
    \begin{array}{ccccc}
      I_s & {\bf 0} & {\bf 0} & {\bf 0} & {\bf 0} \\
      {\bf 0} & J_s-I_s & {\bf 0} & {\bf 0} & {\bf 0} \\
      {\bf 0} & {\bf 0} & 1-\frac{s}{s-1}-\frac{t}{t-1} & {\bf 0} & {\bf 0} \\
      {\bf 0} & {\bf 0} & {\bf 0} & J_t-I_t & {\bf 0} \\
      {\bf 0} & {\bf 0} & {\bf 0} & {\bf 0} & I_t \\
    \end{array}
  \right).$$
Clearly, $P^T(A(H^l)+I)P$ is invertible, then $m_{A(H^l)}(-1)=0$, and thus $m_{L(H)}(1)=0$ by Lemma \ref{eqlemma}.
\hfill$\square$

\section{Proofs of Theorems \ref{mainth0}, \ref{mainth1} and \ref{mainth2}}

\quad Now we are ready to show the proofs of Theorems \ref{mainth0}, \ref{mainth1} and \ref{mainth2}.

\vskip 2mm
\noindent
{\bf Proof of Theorem \ref{mainth0}:} \ If $p(G)=q(G)$, then $G=\overline{G}$, and the conclusion is clear. Next, suppose that $p(G)>q(G)$. Then there exists at least one  quasi-pendant vertex possessing greater than one pendant vertex and let $\{v_1, v_2, \cdots, v_t\}$ be all such quasi-pendant vertices in $G$. We first consider the pendant vertices $\{w_1, w_2, \cdots, w_k\}$ $(k>1)$ adjacent to $v_1$. By using Lemma \ref{mainlemma} to the edge $e_{v_1w_1}$, we have
$$m_{L(G)}(1)=m_{L(P_3)}(1)+m_{L(G-w_1)}(1)=1+m_{L(G-w_1)}(1).$$
Denote by $G_1=G-\{w_1, w_2, \cdots, w_{k-1}\}$. Then recursively using Lemma \ref{mainlemma} to the edges $e_{v_2w_2}$, $e_{v_3w_3}$, $\cdots$, $e_{v_{k-1}w_{k-1}}$, it follows that $m_{L(G)}(1)=k-1+m_{L(G_1)}(1).$ Further, One could consider $v_2$, $v_3$, $\cdots$, $v_t$ one by one analogous with above, then the conclusion is clear. \hfill$\square$

\vskip 2mm
\noindent
{\bf Proof of Theorem \ref{mainth1}:} \ Let $T\in \mathcal{T}(n)$. Since $n\geq 6$ and $p(T)=q(T)$, then the diameter of $T$ is at least 4.  Let $P_k=v_1v_2\cdots v_k$ $(k\geq 5)$ be a diametrical path of $T$. Then $d_{v_2}=d_{v_{k-1}}=2$, and $d_{v_3}>2$ and $d_{v_{k-2}}>2$ as $T$ contains no pendant path $P_3$. If $6\leq n\leq 9$, then there are just 7 trees (see Fig. 5). By using Lemmas \ref{path3}, \ref{mainlemma} and \ref{starlike}, one can easily obtain that $m_{L(T)}(1)=0 \leq \frac{n-6}{4}$ for the trees in Fig. 5, and $m_{L(T)}(1)=0=\frac{n-6}{4}$ if and only if $T$ is isomorphic to $T_1$ of Fig. 5, as required.

\begin{figure}[htbp]
  \centering
  \setlength{\abovecaptionskip}{0cm}
  \setlength{\belowcaptionskip}{0pt}
  \includegraphics[width=4 in]{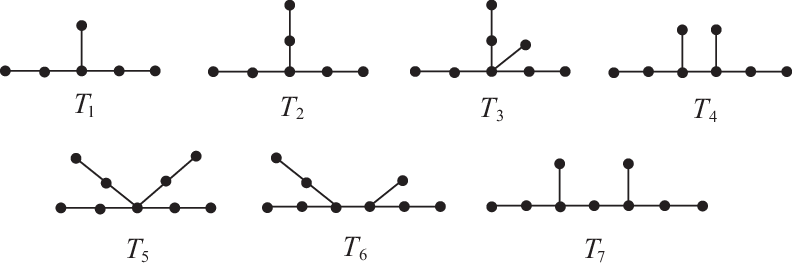}
  \caption{The trees of $\mathcal{T}(n)$ with $6\leq n\leq 9$.}
\end{figure}

In the following, suppose that $n\geq 10$ and the conclusion holds for all trees with order less than $n$. We mainly use induction on the order $n$ to show the remaining proof, which can be divided into two cases.

\vskip 2mm
\noindent
{\bf Case 1.}\ Neither $v_3$ nor $v_{k-2}$ in the diametrical path $P_k$ is a quasi-pendant vertex.

For this case, apart from the diametrical path $P_k$, there can only be some pendant path $P_2$ adjacent to $v_3$ and $v_{k-2}$. As $d_{v_3}>2$ (resp., $d_{v_{k-2}}>2$), then $v_3$ (resp., $v_{k-2}$) possesses at least one pendant $P_2$. If $diam(T)=4$, then $T$ is a star-like tree and $m_{L(T)}(1)=0$ from Lemma \ref{starlike}. Now assume that $diam(T)\geq 5$. Denote by $T_{v_3}$ (resp., $T_{v_4}$) the component of $T-e_{v_3,v_4}$ (removing the edge $e_{v_3,v_4}$ from $T$), which contains the vertex $v_3$ (resp., $v_4$). Clearly, $T_{v_3}$ is a star-like tree, and then from Lemma \ref{starlike}
\begin{equation}\label{e2}
  m_{L(T_{v_3})}(1)=0.
\end{equation}

In what follows, we mainly consider the structure of $T_{v_4}$: (i)\ $T_{v_4}$ contains no pendant $P_3$ and $p(T_{v_4})=q(T_{v_4})$ (i.e., $T_{v_4}\in \mathcal{T}(n)$); (ii)\ $T_{v_4}$ contains pendant $P_3$; (iii)\ $p(T_{v_4})>q(T_{v_4})$.  We shall subsequently show the proofs for the three subcases.

${\bf (i)}$ \ If $T_{v_4}\in \mathcal{T}(n)$, then $|T_{v_4}|\geq 6$. Applying induction hypothesis to $T_{v_4}$, we get
$$m_{L(T_{v_4})}(1)\leq \frac{|T_{v_4}|-6}{4}.$$
Note that $|T_{v_4}|=n-|T_{v_3}|\leq n-5$, then it follows from Lemma \ref{boundlemma} that
\begin{equation*}
\begin{array}{rcl}
m_{L(T)}(1) & \leq & 1+m_{L(T_{v_3})}(1)+m_{L(T_{v_4})}(1)=1+m_{L(T_{v_4})}(1)\\
&\leq & \frac{|T_{v_4}|-2}{4} \leq \frac{n-7}{4}< \frac{n-6}{4}.
\end{array}
\end{equation*}

${\bf (ii)}$ \ Suppose that $T_{v_4}$ contains pendant path $P_3$, then $T$ is isomorphic to one of the trees $\{H_1, H_2, H_3\}$ of Fig. 6 with the depicted local structure.

Firstly, if $T$ is isomorphic to $H_1$ of Fig. 6, then denote by $\overline{H}_1$ the tree obtained from $H_1$ by deleting $\{v_4, v_5, v_6\}$ and joining $v_3$ and $v_7$ with an edge, and thus $m_{L(T)}(1)=m_{L(\overline{H}_1)}(1)$ from Corollary \ref{innerpathcorol}. Using induction hypothesis to $\overline{H}_1$, we have
$$m_{L(T)}(1)=m_{L(\overline{H}_1)}(1)\leq \frac{|\overline{H}_1|-6}{4}=\frac{n-9}{4}< \frac{n-6}{4}.$$

Secondly, if $T$ is isomorphic to $H_2$ of Fig. 6, then applying Lemma \ref{mainlemma} to $u_1, v_4$ of $H_2$, it is obtained that
\begin{equation}\label{e3}
m_{L(T)}(1)=m_{L(T_{v_3}')}(1)+m_{L(T_{v_4})}(1),
\end{equation}
where $T_{v_3}'$ is the reduction component containing $T_{v_3}$.  Obviously, $T_{v_3}'$ is also a star-like tree and then $m_{L(T_{v_3}')}(1)=0$ from Lemma \ref{starlike}. Recalling that $v_{k-2}$ in the diametrical path $P_k$ is not a quasi-pendant vertex, then $v_4\neq v_{k-2}$ and $|T_{v_4}|\geq 6$. Since $T_{v_4}$ contains a pendant $P_3$ and still satisfies $p(T_{v_4})=q(T_{v_4})$,  then Lemma \ref{previousth} tells us that
$$m_{L(T_{v_4})}(1)\leq \frac{|T_{v_4}|-2}{4}.$$
Hence, by equation (\ref{e3}),
$$m_{L(T)}(1)=m_{L(T_{v_4})}(1)\leq \frac{|T_{v_4}|-2}{4} \leq \frac{n-5-2}{4}< \frac{n-6}{4}.$$

Thirdly, suppose that $T$ is isomorphic to $H_3$ of Fig. 6. If $|T_{v_4}|=5$, then $T_{v_4}=P_5$ is a star-like tree and $T$ can be viewed as the tree $H$ of Lemma \ref{starlike}, and thus $m_{L(T)}(1)=0$. If $|T_{v_4}|\neq 5$, then $|T_{v_4}|\geq 8$. Let  $|T_{v_4}|= 8$, then we say that $T$ is isomorphic to $H_4$ of Fig. 6 and Lemma \ref{path3} indicates that $m_{L(T_{v_4})}(1)=0$. By Lemma \ref{boundlemma} and (\ref{e2}), we further have
$$m_{L(T)}(1) \leq 1+m_{L(T_{v_3})}(1)+m_{L(T_{v_4})}(1)=1<\frac{n-6}{4}.$$
Now suppose that $|T_{v_4}|> 8$. Denote $\overline{T}_{v_4}=T_{v_4}-\{u_1,u_2,v_4\}$ in $H_3$, then from Lemma \ref{path3}
$$m_{L(T_{v_4})}(1)=m_{L(\overline{T}_{v_4})}(1).$$
If $p(\overline{T}_{v_4})=q(\overline{T}_{v_4})$, then from Lemma \ref{previousth} (noting that $\overline{T}_{v_4}$ cannot satisfy the conditions of extremal trees, as $d_{v_{k-2}}>2$ and $v_{k-2}$ is not a quasi-pendant vertex)
$$m_{L(\overline{T}_{v_4})}(1)< \frac{|\overline{T}_{v_4}|-2}{4},$$
which implies that
$$m_{L(T)}(1) \leq 1+m_{L(T_{v_3})}(1)+m_{L(T_{v_4})}(1)< 1+\frac{|\overline{T}_{v_4}|-2}{4}\leq  1+\frac{(n-8)-2}{4}=\frac{n-6}{4}.$$
If $p(\overline{T}_{v_4})>q(\overline{T}_{v_4})$, then $T$ must have the local structure as depicted in $H_5$ of Fig. 6. We turn to considering the two components $T_{v_5}$ and $T_{v_6}$ of $T-e_{v_5,v_6}$, which contains $v_5$ and $v_6$, respectively. Since $v_6$ is a quasi-pendant vertex, then $v_6$ cannot be the vertex $v_{k-2}$, and so $|T_{v_6}|\geq 7$. It follows from Lemma \ref{previousth} and $p(T_{v_6})=q(T_{v_6})$ that
\begin{equation}\label{e4}
m_{L(T_{v_6})}(1)\leq \frac{|T_{v_6}|-2}{4}.
\end{equation}
Applying Lemma \ref{mainlemma} to the vertices $u_3$ and $v_6$ in $H_5$, we get
\begin{equation}\label{e5}
m_{L(T)}(1)=m_{L(T_{v_5}')}(1)+m_{L(T_{v_6})}(1),
\end{equation}
where $T_{v_5}'$ is the reduction component containing $T_{v_5}$. Using Lemmas \ref{path3} and \ref{starlike}, we see that
\begin{equation}\label{e6}
m_{L(T_{v_5}')}(1)=0.
\end{equation}
Combining (\ref{e4}), (\ref{e5}) and (\ref{e6}), we derive that
$$m_{L(T)}(1)\leq \frac{|T_{v_6}|-2}{4}\leq \frac{(n-9)-2}{4}<\frac{n-6}{4}.$$

\begin{figure}[htbp]
  \centering
  \setlength{\abovecaptionskip}{0cm}
  \setlength{\belowcaptionskip}{0pt}
  \includegraphics[width=5.7 in]{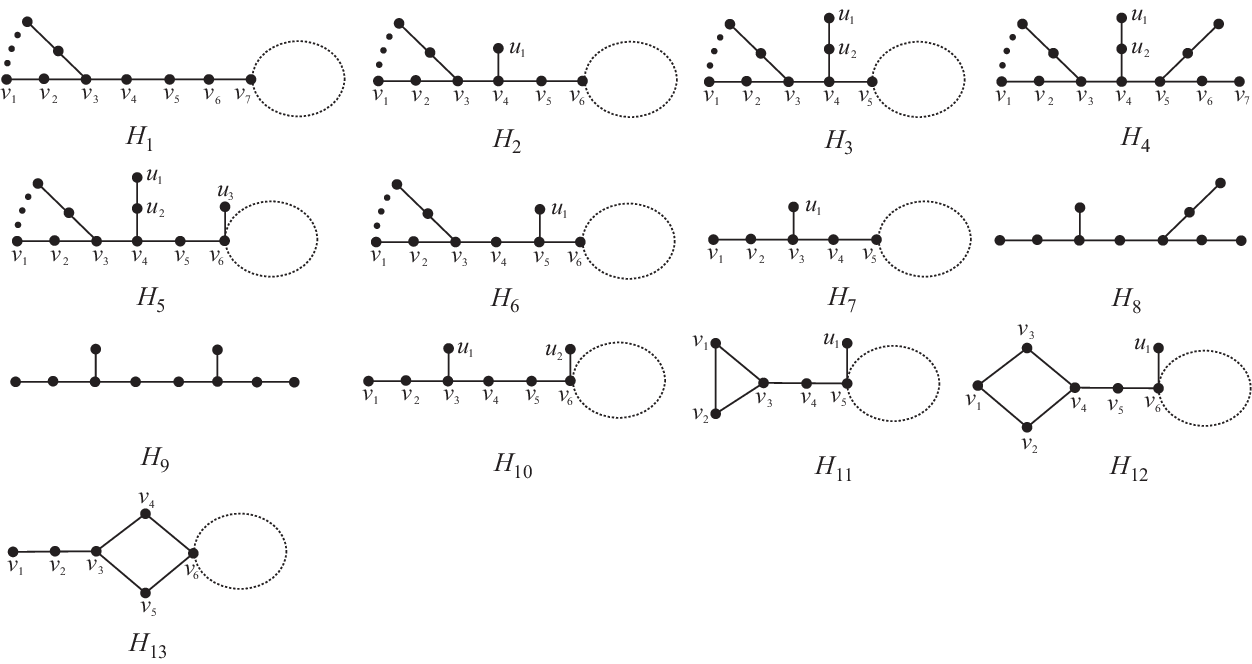}
  \caption{The graphs $H_1, H_2,\cdots H_{13}$ (the dashed circle means the remaining part).}
\end{figure}

${\bf (iii)}$ \ Suppose that $p(T_{v_4})>q(T_{v_4})$, then $T$ has the local structure as depicted in $H_6$ of Fig. 6. Denote by  $T_{v_4}$ and $T_{v_5}$ the two components of $T-e_{v_4,v_5}$, which contains $v_4$ and $v_5$, respectively. Do reduction operation to $H_6$ based on $u_1$ and $v_5$, and let $T_{v_4}'$ be the reduction component containing $T_{v_4}$. By Lemmas \ref{path3} and \ref{starlike},
$$m_{L(T_{v_4}')}(1)=0.$$
Applying Lemma \ref{previousth}, we get
$$m_{L(T_{v_5})}(1)\leq \frac{|T_{v_5}|-2}{4}.$$
Then it follows from Lemma \ref{mainlemma} that
\begin{equation*}
m_{L(T)}(1)=m_{L(T_{v_4}')}(1)+m_{L(T_{v_5})}(1)\leq \frac{|T_{v_5}|-2}{4}\leq \frac{(n-6)-2}{4}<\frac{n-6}{4}.
\end{equation*}

\vskip 2mm
\noindent
{\bf Case 2.} \ Either $v_3$ or $v_{k-2}$ in the diametrical path $P_k$ is a quasi-pendant vertex.

Without loss of generality, suppose that $v_3$ is a quasi-pendant vertex. If $d_{v_3}\geq 4$, then $v_3$ possesses at least two  pendant path $P_2$. Recalling that $n\geq 10$, one can easily obtain by Lemma \ref{mainlemma} that
\begin{equation}\label{e7}
m_{L(T)}(1)=m_{L(P_4)}(1)+m_{L(T-P_2)}(1)=m_{L(T-P_2)}(1).
\end{equation}
Clearly, $T-P_2\in \mathcal{T}(n)$ still, then from induction hypothesis to $T-P_2$,
\begin{equation}\label{e8}
m_{L(T-P_2)}(1)\leq \frac{|T-P_2|-6}{4}=\frac{n-8}{4}.
\end{equation}
Combining (\ref{e7}) and (\ref{e8}), we get
\begin{equation}\label{e9}
m_{L(T)}(1)\leq \frac{n-8}{4}<\frac{n-6}{4}.
\end{equation}

Now suppose that $d_{v_3}=3$. For brevity, denote $T_{v_3}=T-\{v_1,v_2\}$. Similar as (\ref{e7}), we obtain
\begin{equation}\label{e10}
m_{L(T)}(1)=m_{L(T_{v_3})}(1).
\end{equation}
We further consider the structure of $T_{v_3}$. If $T_{v_3}$ contains no pendant $P_3$, then $T_{v_3}\in \mathcal{T}(n)$ and parallel with the above discussion we also have the inequality (\ref{e9}). If $T_{v_3}$ contains a pendant $P_3$, then $T$ has the local structure as depicted in $H_7$ of Fig. 6. If $|T|=10$, then $T$ is isomorphic to $H_8$ or $H_9$ of Fig. 6. From Lemmas \ref{mainlemma} and \ref{path3}, we can easily obtain that $m_{L(H_8)}(1)=0$ and $m_{L(H_9)}(1)=1=\frac{|H_9|-6}{4}$. Next, assume that $|T|\geq 11$. Denote $\overline{T}_{v_3}=T_{v_3}-\{u_1, v_3, v_4\}$ in $H_7$ and then $|\overline{T}_{v_3}|\geq 6$. If $p(\overline{T}_{v_3})=q(\overline{T}_{v_3})$, then from Lemmas \ref{path3} and \ref{previousth}
$$m_{L(T_{v_3})}(1)=m_{L(\overline{T}_{v_3})}(1)\leq \frac{|\overline{T}_{v_3}|-2}{4}= \frac{(n-5)-2}{4}<\frac{n-6}{4},$$
which, together with (\ref{e10}), implies that
$$m_{L(T)}(1)<\frac{n-6}{4}.$$
If $p(\overline{T}_{v_3})>q(\overline{T}_{v_3})$, then $T$ is isomorphic to $H_{10}$ of Fig. 6.  Denote $T_{v_4}=T-\{v_1,v_2,v_3,u_1\}$ in $H_{10}$, and $T_{v_4}\in \mathcal{T}(n)$ obviously. As a result, from Lemma \ref{extremallemma} and the induction hypothesis to $T_{v_4}$, we have
\begin{equation}\label{e11}
m_{L(T)}(1)=1+m_{L(T_{v_4})}(1)\leq 1+\frac{|T_{v_4}|-6}{4}=1+\frac{(n-4)-6}{4}=\frac{n-6}{4}.
\end{equation}
Therefore, we conclude that $m_{L(T)}(1)\leq \frac{n-6}{4}$ when  $T\in \mathcal{T}(n)$.

At last, we show the sufficiency and necessity of the equality. For the sufficiency part, if $T$ is isomorphic to the tree of Fig. 1, then Lemma \ref{extremallemma} tells us that $m_{L(T)}(1)=\frac{n-6}{4}.$
For the necessity part, suppose that $m_{L(T)}(1)=\frac{n-6}{4}$ with $T\in \mathcal{T}(n)$. From the above discussions, if $|T|=6$ or $|T|=10$, then $T$ is isomorphic to $T_1$ of Fig. 5 or $H_9$ of Fig. 6, and the conclusion holds. If $|T|> 10$, then there can only be the inequality (\ref{e11}) turning to equality, that is,
$$m_{L(T)}(1)=1+m_{L(T_{v_4})}(1)=1+\frac{|T_{v_4}|-6}{4}.$$
The induction hypothesis to $T_{v_4}$ yields $T$ has the structure depicted as the tree in Fig. 1.  The proof of Theorem \ref{mainth1} is completed.
\hfill$\square$

\vskip 2mm
At last, we prove Theorem \ref{mainth2}.

\vskip 2mm
\noindent
{\bf Proof of Theorem \ref{mainth2}:} \ Let $G\in \mathcal{G}(n)$ be a unicyclic graphs on $n(\geq 10)$ vertices. For brevity, we call the vertices on the cycle of $G$ the cycle-vertices. We divide the proof into two cases.

\vskip 2mm
\noindent
{\bf Case 1.} \ There exists a cycle-vertex, say $v$, in $G$ which is a quasi-pendant vertex.

Let $w$ be another cycle-vertex and $w\thicksim v$ in $G$. Applying Lemma \ref{mainlemma} to the edge $e_{vw}$, then
$m_{L(G)}(1)=m_{L(\Gamma)}(1),$
where $\Gamma$ is the graph obtained as in Lemma \ref{mainlemma}. Clearly, $\Gamma$ is a tree with order $|\Gamma|=n+2$ and $p(\Gamma)=q(\Gamma)$. Then by Lemma \ref{previousth},
$$m_{L(\Gamma)}(1)\leq \frac{|\Gamma|-2}{4}=\frac{n}{4}.$$
Hence, it follows that
\begin{equation}\label{e12}
m_{L(G)}(1)=m_{L(\Gamma)}(1)\leq \frac{|\Gamma|-2}{4}=\frac{n}{4}.
\end{equation}

\vskip 2mm
\noindent
{\bf Case 2.} \ Each cycle-vertex in $G$ is not a quasi-pendant vertex.

For this case, if $G$ is a cycle $C_n$, then it is known that the Laplacian eigenvalues of $C_n$ are $2-2cos\frac{2\pi j}{n}$ $(j=0, 1, \cdots, n-1)$. Thus $m_{L(C_n)}(1)=0$ or $2$, and $m_{L(C_n)}(1)=2$ if $6$ divides $n$. Recalling that $n\geq 9$, then we say $m_{L(G)}(1)< \frac{n}{4}$.

Next, suppose that $G\neq C_n$ and $v$ is a cycle-vertex with $d_v\geq 3$ in $G$. Let $w$ be also a cycle-vertex of $G$ and $w\thicksim v$. Denote by $T=G-e_{vw}$ the tree obtained from $G$ by deleting the edge $e_{vw}$. It is not hard to see that $p(T)=q(T)$.
Moreover, if $T$ contains no pendant path $P_3$, then by Lemma \ref{boundlemma} and Theorem \ref{mainth1} we get
\begin{equation}\label{e13}
m_{L(G)}(1)\leq 1+m_{L(T)}(1)\leq 1+\frac{n-6}{4}=\frac{n-2}{4}<\frac{n}{4}.
\end{equation}
Now suppose that $T$ contains at least one pendant path $P_3$. 
If $p(T-P_3)=q(T-P_3)$ holds, then it follows from Lemmas \ref{boundlemma}, \ref{path3} and \ref{previousth} that
$$m_{L(G)}(1)\leq 1+m_{L(T)}(1)=1+m_{L(T-P_3)}(1)\leq 1+\frac{(n-3)-2}{4}=\frac{n-1}{4}<\frac{n}{4}.$$
If $p(T-P_3)>q(T-P_3)$, then $G$ has the local structure as depicted in $H_{11}$, $H_{12}$ or $H_{13}$ of Fig. 6. If $G$ is isomorphic to $H_{11}$, then let $\Gamma$ be the graph obtained from $H_{11}-e_{v_4v_5}$ and a disjoint path $P_2$ by
joining $v_4$ to a vertex of $P_2$. Clearly, $\Gamma$ consists of two disjoint components: a unicyclic graph of order $6$ (denoted by $G_1$) and a reduced tree of order $n-4$ (denoted by $T_1$). It is easy to see that $m_{L(G_1)}(1)=0$, and then by Lemmas \ref{mainlemma} and \ref{previousth}
$$m_{L(G)}(1)=m_{L(\Gamma)}(1)=m_{L(T_1)}(1)\leq \frac{|T_1|-2}{4}=\frac{n-6}{4}<\frac{n}{4}.$$
If $G$ is isomorphic to $H_{12}$, then we could delete the edge $e_{v_1v_3}$ from $H_{12}$ and the resultant graph is still reduced and without pendant path $P_3$. Thus we could also obtain the inequality (\ref{e13}). Similarly, if $G$ is isomorphic to $H_{13}$, then we could also get the inequality (\ref{e13}) by deleting $e_{v_4v_6}$ from $H_{13}$.

Combining the above two cases, the inequality in Theorem \ref{mainth2} is proved. At last, we show the sufficiency and necessity for the equality. If $G$ is isomorphic to the graph in Fig. 2, then we easily obtain $m_{L(G)}(1)=\frac{n}{4}$ by applying Corollary \ref{corol} recursively. For the necessity part, let
$m_{L(G)}(1)=\frac{n}{4}$, then $G$ belongs to Case 1,  i.e., the inequality (\ref{e12}) turns to equality. Thus  the tree $\Gamma$ in Case 1 satisfies the two assertions of Lemma \ref{previousth}. Let $v$ be a quasi-pendant cycle-vertex of $G$ as assumed in Case 1. Since $\Gamma$ satisfies the assertion (i) of Lemma \ref{previousth}, then each cycle-vertex, distinct with $v$, of degree at least 3 is a quasi-pendant vertex.
Consider the degree of the vertex $v$ in Case 1. We say $d_v=3$ in $G$, that is, there is only a pendant vertex adjacent to $v$ apart from two cycle-vertices in $G$. Otherwise, assume that $d_v\geq 4$, then there exists a component, adjacent to $v$, containing pendant path $P_4$ (noting that $\Gamma$ satisfies the assertion (ii) of Lemma \ref{previousth}), a contradiction. Hence we could say that the degree of  each quasi-pendant cycle-vertex of $G$ is 3. From the above, we could claim that $G$ is obtained from a cycle by adding some pendant vertices to some cycle-vertices such that each cycle-vertex is of degree 2 or 3. Applying the assertion (ii) of Lemma \ref{previousth} to $\Gamma$ again, we get that any two nearest cycle-vertices in distance are connected by a path $P_2$ in $G$. Consequently, $G$ is isomorphic to the graph in Fig. 2. The proof of Theorem \ref{mainth2} is completed.
\hfill$\square$

\end{document}